\documentstyle{amsppt} 

\def \al{\alpha}
\def \bt{\beta}

\def \KK{\Bbb {K}}

\def \NN{\Bbb {N}}
\def \ZZ{\Bbb {Z}}

\def \sd #1#2{{#1}^{\scriptstyle{(#2)}}}

\def \pd #1#2{#1^{(#2)}}
\def \gna #1#2{\Gamma_{#1}(#2)^{\scriptstyle {ab}}}
\def\ns#1{\vskip10pt\np{\bf #1}\qquad}


\hoffset-8pt
\voffset 20pt

\font\fintr=cmr9

\font\smc=cmtt10

\magnification=\magstep1

\hsize=6.5truein

\vsize=8.5truein

\NoBlackBoxes

\def\np{\par\noindent}

\def\endemo{\hfill\qed\enddemo}

\document
\topmatter
\title Homogeneous Multiplicative Polynomial Laws are Determinants\endtitle
\author F. Vaccarino \endauthor
\address Dipartimento di  Matematica - Politecnico di Torino
- Corso Duca degli Abruzzi 24 - 10129 - Torino - Italy
\quad e-mail: vaccarino \@ syzygie.it
\endaddress
\endtopmatter

\vskip20pt
{\fintr \smallskip\leftskip=25pt \rightskip=25pt
\baselineskip8pt  \centerline{\smc{Abstract}}
Let $R$ be a ring and let $B$ be a commutative ring. Let $p:R @>>> B$ be a 
homogeneous multiplicative polynomial law of degree $n$. 
The main result of 
this paper is to show that $p$ is essentially a 
determinant, in the sense that $p$ is obtained from a 
determinant by left and right 
composition with ring homomorphisms. 
This is achieved using results on the invariants of 
matrices in positive characteristic due to S.\,Donkin {\cite{3}}, 
{\cite{4}},  A.\,Zubkov {\cite{9}}, {\cite{10}} and the author 
{\cite{7}}. \par}

\vskip20pt
\head Introduction \endhead
\vskip10pt
Let $R$ be a ring and $B$ a commutative ring. Let $p:R @>>> B$ be a 
homogeneous multiplicative polynomial law of degree $n$. The 
definition of polynomial law will be recalled below and the reader 
can think about 
this as a polynomial map between rings that preserves the product 
and the 
identity. We want to understand the links between $p$ and the usual 
determinant of $n\times n$ matrices, which is itself a homogeneous 
multiplicative 
polynomial law of degree $n$. 

In order to state the main result of this paper we have to introduce 
some objects.
\medskip

Let $S$ be a set and let $F_S:=\ZZ \langle x_{s}\rangle_{s\in S}$ be 
the free ring
on it. Let $\ZZ[x_{i,j}^{s}]_{1\leq i,j \leq n, \, s\in S}$ be the
ring of polynomials in the variables $x_{i,j}^{s}$ with $1\leq i,j 
\leq n, \, s\in
S$ and let $M_{n}(\ZZ[x_{i,j}^{s}])$ denote the ring of $n\times
n$ matrices with entries in it.
  
Let $j_{n}: F_{S} @>>> M_{n}(\ZZ[x_{i,j}^{s}])$
 be given by 
$j_n(x_s)=\sum_{i,j=1}^n x_{i,j}^s e_{i,j}$, where 
$(e_{i,j})_{h,k}=\delta_{i,h}\delta_{j,k}$ for $i,j,h,k=1,\dots, n$. 

Let $E_S(n)$ be the subring of $\ZZ[x_{i,j}^{s}]_{1\leq i,j \leq n,
 \, s\in S}$
generated by the coefficients of the characteristic polynomials of the 
elements of $j_{n}(F_{S})$. Then $\det \cdot j_n$ maps $F_S$ to $E_S(n)$.

Let $G=GL_{n}(\ZZ)$ act by simultaneous conjugation on the direct sum 
of ${\scriptstyle{\#}} S$ copies 
of $M_{n}(\ZZ)$. Since $\ZZ[x_{i,j}^{s}]_{1\leq i,j \leq n, \, s\in
S}$
 is the
symmetric algebra of the direct sum 
of ${\scriptstyle{\#}} S$ copies 
of $M_{n}(\ZZ)$, then $G$ acts on it as an automorphism 
group. We denote by $\ZZ[x_{i,j}^{s}]^{\scriptscriptstyle{G}}$ the 
ring of invariants for 
this action. Then, clearly $\ZZ[x_{i,j}^{s}]^{\scriptscriptstyle{G}}$
contains $E_S(n)$. In fact, by works of  S.\, Donkin and A.\, 
Zubkov ({\cite{3}}, {\cite{4}}, {\cite{9}} and {\cite{10}}), we have 
$$\ZZ[x_{i,j}^{s}]^{\scriptscriptstyle{G}}= E_S(n).$$ 

\bigskip
We are now able to state the main result of this paper. 
\smallskip
\proclaim{Theorem}
Let $B$ be a commutative ring and let $p:F_{S} @>>> B$ be a 
homogeneous multiplicative polynomial law of degree $n$. 
There is a unique ring homomorphism 
$\phi:  \ZZ[x_{i,j}^{s}]^{\scriptscriptstyle{G}}@>>> B$
such that the following diagram commutes:

$$
\CD  F_S @>p>> B \\
@V j_n VV        @AA\phi A \\
j_n(F_S) @> \det >> \ZZ[x_{i,j}^{s}]^{\scriptscriptstyle{G}},
\endCD
$$
i.e. $p$ is a determinant composed with ring homomorphisms.
\endproclaim
\medskip

Let $R$ be a ring and let  
$$0 @>>> K @>>> F_S @>\pi >> R @>>>0$$ 
be any presentation of $R$ by means of generators and relations.

Let $B$ be a commutative ring and let $p:R @>>> B$ be a 
homogeneous multiplicative polynomial law of degree $n$.
Now $\pi \cdot p :F_S @>>> B$ is a homogenous multiplicative
polynomial law of degree $n$, by the previous theorem there is a unique 
$\rho :  \ZZ[x_{i,j}^{s}]^{\scriptscriptstyle{G}} @>>> B$ such that
 
$$p(r)=\big (p \cdot \pi \big )(f)=\big (\rho \cdot \det \cdot 
j_{n} \big )(f),$$ 
for all $r\in R$, where $f\in F_{S}$ is such that $\pi(f)=r$, 
i.e. $p$ is a determinant composed with ring homomorphisms.

\bigskip
The paper is divided into three sections: in the first one we recall 
some generalities on polynomial laws and divided powers algebras.

In the second section we recast some results due to S.\, Donkin, A.\, 
Zubkov and the author on the invariants of matrices 
in positive characteristic.  

In the third section we give the main result.

I would like to thank C.Procesi and M.Brion for useful discussions.

\vskip20pt

\head \S 0 Conventions and Notations \endhead

\vskip10pt

\noindent
Except otherwise stated all the rings (algebras over a
commutative ring) should  be understood associative and
with multiplicative identity.

\noindent
We denote by $\NN$ the set of non-negative integers and by $\ZZ$ the 
ring of integers. 

\noindent
Let $R$ be a commutative ring and let $n$ be any positive integer: 
we  denote by $M_{n}(R)$ the ring of $n\times n$ matrices
with entries  in $R$. 

\noindent
Let $A$ be a set, we denote by ${\scriptstyle{\#}} A$ its
cardinality. 

\noindent
For $A$ a set and any additive monoid $M$, we denote by
$M^{\scriptstyle{(A)}}$ the  set of functions $f:A@>>> M$
with finite support.  

\noindent
Let $\al\in \sd M A$, we denote by $\mid \al \mid$ the
(finite)  sum $\sum_{a\in A} \al(a)$.

\noindent
Let $S$ be a set. We put 

$F_S:=\ZZ\langle x_s \rangle_{s\in S}$ for the free ring over 
$S$;

$F_{S}^{+}$ for the augmentation ideal of $F_{S}$ (non-unital 
subring);

$P_{S}:=\ZZ [x_s ]_{s\in S}$ for the free commutative ring over $S$.

\vskip20pt
\head \S 1 Polynomial Laws and Divided Powers \endhead
\vskip10pt
\ns{1.1}
Let $\KK$ be a commutative ring. 
Let us recall the definition of a kind of map between $\KK$-modules
that generalizes the concept of polynomial map between free 
$\KK$-modules (see {\cite{2}} or {\cite{5}}).

\definition{Definition 1.1.1 }
Let $A$ and $B$ be two $\KK$-modules. A {\it{polynomial law }}
$\varphi$ from $A$ to $B$ is a family of mappings 
$\varphi_{_{L}}:L\otimes_{\KK} A \longrightarrow L\otimes_{\KK} B$, 
with $L$ varying in the family of  commutative $\KK$-algebras, such 
that the following diagram commutes:

$$
\CD L\otimes_{\KK}A @> \varphi_{_{L}} >> L\otimes_{\KK} B \\
    @V f\otimes 1_{_{A}}VV           @VV f \otimes 1_{_{B}}V \\
    M \otimes_{\KK} A @> \varphi_{_{M}}>> M\otimes_{\KK}B,
    \endCD
    $$ 
for all $L$, $M$ commutative $\KK$-algebras and all  
homomorphisms of $\KK$-algebras $f:L \longrightarrow M$.
\enddefinition
\smallskip
\definition{Definition 1.1.2}
Let $n\in \NN$, if $\varphi_{_{L}}(au)=a^n\varphi_{_{L}}(u)$, 
for all $a\in L$, $u\in L\otimes_{\KK} 
A$ and all commutative $\KK$-algebras $L$, then $\varphi$ will 
be said {\it{homogeneous of degree $n$}}.
\enddefinition
\smallskip

\definition{Definition 1.1.3}
If $A$ and $B$ are two $\KK$-algebras  and 
$$\cases \varphi_{_{L}}(xy)&=\varphi_{_{L}}(x)\varphi_{_{L}}(y)\\
         \varphi_{_{L}}(1_{_{L\otimes A}})&=1_{_{L\otimes B}},
         \endcases
         $$
for all commutative $\KK$-algebras $L$ and for all 
$x,y\in L\otimes A$, then $\varphi$ is called 
{\it{multiplicative}}.
\enddefinition
\smallskip

Let $A$ and $B$ be two $\KK$-modules and $\varphi:A@>>>B$ be a 
polynomial law. 
We recall the following result on polynomial laws, which is a 
restatement of Th\'eor\`eme I.1 of {\cite{5}}.
\proclaim{Proposition 1.1.4} Let $S$ be a set.
\roster
\item Let $L=\KK\otimes P_{S}$ and let $\{a_{s}\, :\,s\in S\}\subset A$
be such that $a_{s}=0$ except for a finite number of $s\in S$, then
there exist $\varphi_{\xi}((a_{s}))\in B$, with $\xi \in \sd {\NN} S$,
such that:
$$\varphi_{_{L}}(\sum_{s\in S} x_s\otimes a_{s})=\sum_{\xi \in
\NN^{(S)}}  x^{\xi}\otimes \varphi_{\xi}((a_{s})),$$ where 
$x^{\xi}:=\prod_{s\in S} 
x_s^{\xi_s}$.
\item Let $R$ be any commutative $\KK$-algebra and let 
$(r_s)_{s\in S}\subset R$, then:
$$\varphi_{_{R}}(\sum_{s\in S} r_s\otimes a_{s})=\sum_{\xi \in
\NN^{(S)}}  r^{\xi}\otimes \varphi_{\xi}((a_{s})),$$ where 
$r^{\xi}:=\prod_{s\in S} 
r_s^{\xi_s}$.
\item If $\varphi$ is homogeneous of degree $n$, then in the previous 
sum one has $\varphi_{\xi}((a_{s}))=0$ if $\mid \xi \mid$ is 
different from $n$. That is:
$$\varphi_{_{R}}(\sum_{a\in A} r_a\otimes a)=\sum_{\xi \in 
\NN^{(A)},\,\mid \xi \mid=n} r^{\xi}\otimes \varphi_{\xi}((a)).$$
In particular, if $\varphi$ is homogeneous of degree $0$ or $1$, then 
it is constant or linear, respectively.
\endroster
\endproclaim
\smallskip
Let $S$ be a set, Proposition 1.1.4 means that a polynomial law 
$\varphi:A@>>>B$ is 
completely determined by its coefficients $\varphi_{\xi}((a_{s}))$, with
$(a_s)_{s\in S} \in \sd A S$.  

\remark{Remark 1.1.5}
If $A$ is a free $\KK$-module and $\{a_{t}\, :\, t\in T\}$ is a basis 
of $A$, then $\varphi$ is completely determined by its coefficients 
$\varphi_{\xi}((a_{t}))$, with $\xi \in \NN^{(T)}$.
If also $B$ is a free $\KK$-module with basis $\{b_{u}\, :\, u\in 
U\}$, then $\varphi_{\xi}((a_{t}))=\sum_{u\in U}\lambda_{u}(\xi)b_{u}$. 
Let $a=\sum_{t\in T}\mu_{t}a_{t}\in A$. Since only a finite number of 
$\mu_{t}$ and $\lambda_{u}(\xi)$ are 
different from zero, the following makes sense:
$$\varphi(a)=\varphi(\sum_{t\in T}\mu_{t}a_{t})=
\sum_{\xi\in \NN^{(T)}} \mu^{\xi}\varphi_{\xi}((a_{t}))=
\sum_{\xi\in \NN^{(T)}} \mu^{\xi}(\sum_{u\in 
U}\lambda_{u}(\xi)b_{u})=\sum_{u\in U}(\sum_{\xi\in 
\NN^{(T)}}\lambda_{u}(\xi) \mu^{\xi})b_{u}.$$
Hence, if both $A$ and $B$ are free $\KK$-modules, a polynomial law 
between them is simply a polynomial map.
\endremark
\vskip10pt
\ns{1.2}
Let $\KK$ be any commutative ring with identity. 
For a $\KK$-module $M$ let $\Gamma(M)$ denote its divided
powers algebra (see {\cite{2}}, {\cite{5}}). 
This is a unital commutative $\KK$-algebra, with generators $\pd m k $, 
with $m\in M$, $k \in \ZZ$ and relations, for all $m,n\in M$: 
$$
\align 
\pd m i &= 0, \qquad \forall i<0; \tag i \\
\pd m 0 &=1_{\scriptstyle{\Bbb K}}, \qquad \forall m\in M; \tag ii \\
\pd {(rm)} i &= r^i\pd m i, \qquad \forall r\in R, \forall i\in \NN; 
\tag iii \\
\pd {(m+n)} k &= \sum_{i+j=k}\pd m i \pd n j , \qquad \forall k\in
\NN; 
\tag iv \\
\pd m i \pd m j &= {i+j\choose i}\pd m {i+j} , 
\qquad \forall i,j\in \NN. \tag v 
\endalign
$$

The $\KK$-module $\Gamma(M)$ is generated by products 
(over arbitrary index sets $I$) $\prod_{i\in I} \pd {x_i} {\al_i}$
of the above generators, 
it is clear that $\prod_{i\in I} \pd {x_i} {\al_i}=0$ if $\al_i<0$ 
for some $i\in I$.  
The divided powers algebra $\Gamma(M)$ is a $\NN$-graded algebra
with homogeneous components $\Gamma_k:=\Gamma_k(M)$, ($k\in \NN$), the
submodule generated by 
$\{\prod_{i\in I} \pd {x_i} {\al_i}\; :\;\vert \al\vert=k\}$. Note
that $\Gamma_0\cong \KK$ and $\Gamma_1\cong M$. $\Gamma$ is a functor from
$\KK$-modules to commutative unital graded $\KK$-algebras.

Indeed for any morphism of $\KK$-modules $f:M @>>> N$ there 
exists a unique morphism
of graded  $\KK$-algebras $\Gamma(f):\Gamma(M) @>>> \Gamma(N)$ 
such that
$\Gamma(f)(\pd x n)=\pd {f(x)} n$, for any $x\in M$ and $n\geq 0$. 
From this it follows easily that $\Gamma$ is exact.

Furthermore
$\Gamma(L\otimes_{\KK} M)\cong L\otimes_{\KK} \Gamma(M)$ 
as graded rings by means
of $\pd {(1\otimes x)} n \mapsto 1\otimes \pd x n$.  
Thus the map $\Gamma(f)$ commutes with extensions of scalars. 

If $A$ is a (unital) $\KK$-algebra, 
then $\Gamma_k(A)$ is a (unital) $\KK$-algebra too (see
{\cite{6}}).  
To distinguish the new multiplication on $\Gamma_{k}(A)$ from the one of 
$\Gamma(A)$, we denote it by \lq\lq$\tau_k$\rq\rq. We have:
$$
\align
\prod_{i\in I}\pd {a_i} {\al_i}\;\tau_{k}\; 
\prod_{j\in J}\pd {b_j} {\beta_j}&:=
(\prod_{i\in I}\pd {a_i} {\al_i})\; \tau_k \; 
(\prod_{j\in J}\pd {b_j} {\beta_j})\\
&:=\sum_{(\lambda_{ij})\in M(\al,\beta)}\quad 
\prod_{(i,j)\in I\times J}\pd {(a_i b_j)}
{\lambda_{ij}}, \endalign
$$
where $M(\al,\beta):=\{(\lambda_{ij})\in \NN^{(I\times J)} : \sum_{i\in
I}\lambda_{ij}=\beta_j\;,\forall j\in J\; ;\sum_{j\in
J}\lambda_{ij}=\al_i \;,\forall i\in I\}$ and 
$\prod_{i\in I}\pd {a_i} {\al_i}$, 
$\prod_{j\in J}\pd {b_j} {\beta_j}\in \Gamma_k(A)$.
\smallskip

Let us denote by $\gamma_n:=(\gamma_{n,L})$ the polynomial law given by 
the composition $L\otimes M @>>> \Gamma_n(L\otimes M) @>>> L\otimes 
\Gamma_n(M)$, then $\gamma_n$ is homogeneous of degree $n$.

There is a property proved by Roby in {\cite{6}}, which 
motivates our introduction of divided powers.
\proclaim{Theorem 1.2.1}
Let $A$ and $B$ be two $\KK$ -algebras. 
The set of homogeneous multiplicative polynomial laws of 
degree $n$ from $A$ to $B$ is in bijection with the set of 
all homomorphisms of $\KK$-algebras from $\Gamma_{n}(A)$ to $B$. 
Namely, given any homogeneous multiplicative polynomial law $f:A @>>> 
B$ of degree $n$, 
there exists a unique homomorphism of $\KK$-algebras 
$\phi:\Gamma_n(A) @>>> B$ 
such that $f_{_{L}}=(1_{_{L}}\otimes \phi)\cdot \gamma_{n,L}$, 
for any commutative $\KK$-algebra $L$.
\endproclaim
\vskip10pt
\ns{1.3}
There is another result that we need to recall. 

Let $B$ be a commutative ring and let 
$M_{n}(B)$ be the ring of $n\times n$ 
matrices over $B$. 
Let $b\in M_{n}(B)$ and denote by $e_{i}(b)$ the $i$-th 
coefficient of the characteristic polynomial of $b$, i.e. the trace of 
$\wedge^{i}(b)$. 

Let $R$ be a ring, we denote by $(R)^{\scriptscriptstyle{ab}}$ its 
abelianization, that is, its quotient by the ideal generated by the 
commutators of its elements. 

The following can be found in {\cite{8}}.
\proclaim{Proposition 1.3.1}
The ring $\gna n {M_{n}(B)}$ is isomorphic to $B$. The canonical 
projection $\al_{n}:=\al_n(M_{n}(B)):\Gamma_{n}(M_{n}(B)) @>>> \gna n 
{M_{n}(B)}$ is such that, for all $b\in M_{n}(B)$ and $0\leq i \leq n$, 
$$\al_{n}(\pd 1 {n-i} \pd b i)=e_{i}(b).$$
\endproclaim   

\smallskip
For further readings on these topics we refer to 
{\cite{2}}, {\cite{5}}, {\cite{6}}, {\cite{8}}.

\vskip20pt
\head \S 2 Divided Powers and Invariants of Matrices \endhead
\vskip10pt
\ns{2.1}
Let us introduce the following notation: let $S$ be a set, we put

\noindent
$A_S(n):=\ZZ[x_{ij}^s]_{1\leq i,j \leq n\,, \, s\in S}$, the 
symmetric algebra on the direct sum of 
${\scriptstyle{\#}} S$ copies of $M_{n}(\ZZ)$;

\noindent
$B_S(n):=M_n(A_S(n))$; 

\noindent
$G:=Gl_n(\ZZ);$

\noindent
$C_S(n):=A_S(n)^{\scriptscriptstyle{G}}$, 
the ring of polynomial invariants 
with respect to the action of 
$G$ on $A_S(n)$ induced by simultaneous conjugation on the direct sum 
of ${\scriptstyle{\#}} S$ copies of $M_{n}(\ZZ)$.

\vskip10pt
\ns{2.2}
Let $j_n:F_S @>>> B_S(n)$ be given by $j_n(x_s)=\zeta_s:=
\sum_{i,j=1}^n x_{i,j}^s e_{i,j}$, where 
$(e_{i,j})_{h,k}=\delta_{i,h}\delta_{j,k}$ for $i,j,h,k=1,\dots, n$. 
The $\zeta_s$ are the so-called "generic matrices of order $n$". 
Let $\bt_{n}:\Gamma_{n}(B_{S}(n)) @>>> \gna n {B_S(n)}$ be the canonical 
projection, then $A_{S}(n)\cong \gna n {B_S(n)}$ by Prop.1.3.1. 
We set  $E_{S}(n):=\bt_n(\Gamma_{n}(j_{n}(F_{S})) \hookrightarrow 
A_{S}(n)$, by Prop.1.3.1 it is 
the subring of $A_{S}(n)$ generated by the $e_{i}(j_{n}(f))$, where 
$i\in \NN$ and $f\in F_{S}$. 
By the exactness properties of $\Gamma_{n}$ and 
$(-)^{\scriptscriptstyle{ab}}$ there exists a unique ring 
epimorphism $\pi_n: \Gamma_n(F_S)^{\scriptscriptstyle{ab}}@>>> 
E_{S}(n)$ such that the following diagram commutes:

$$
\CD
\Gamma_n(F_S) @> \Gamma_{n}(j_{n})>> \Gamma_{n}(j_{n}(F_{S}))  \\
@V\al_nVV      @VV \beta_n V \\
\Gamma_n(F_S)^{\scriptscriptstyle{ab}}@> \pi_n >> E_{S}(n).
\endCD
$$
\smallskip
\remark{Remark 2.2.1}
By the previous discussion the polynomial law $j_{n}(F_{S}) @>>> 
E_{S}(n)$ that corresponds to $\bt_{n}$, via theorem 1.2.1, 
is the restriction to $j_{n}(F_{S})$ of the determinant.
\endremark
\smallskip
Recall that $E_S(n)$ is a subring of $C_S(n)$.

The Donkin-Zubkov Theorem on invariant of matrices can 
be stated in the following way , see {\cite{3}}, {\cite{4}}, 
{\cite{9}} and {\cite{10}}. (This Theorem was firstly proved by 
S.\;Donkin and then by A.\;Zubkov in another way).

\proclaim{Theorem (Donkin-Zubkov)}
The ring $C_{S}(n)$ of polynomial invariants of the direct sum 
of ${\scriptstyle{\#}} S$ copies of $n\times n$ matrices is equal to 
$E_{S}(n)$. 
\endproclaim

Then we have a surjection $\pi_n:\gna n {F_S}@>>>E_S(n)=C_S(n)$.

The following result is proved in {\cite{7}}.
\proclaim {Theorem 2.2.2}
The map $\pi_n:\gna n {F_S}@>>>C_S(n)$ is an isomorphism of graded rings.
\endproclaim
\smallskip
The proof of the theorem goes as follows: one sees that the natural 
multidegree of $F_S$ induces a multidegree on 
$\Gamma_{n}(F_{s})^{\scriptscriptstyle{ab}}$ and the following is showed
to be an homomorphism of graded
rings for all $n$: 
$$\rho_n:\cases\Gamma_{n}(F_{s})^{\scriptscriptstyle{ab}}
@>>>\Gamma_{n-1}(F_{s})^{\scriptscriptstyle{ab}}\\
\pd 1 {n-\vert \al \vert}\prod_{i\in I}\pd {a_i} {\al_i} \mapsto
\pd 1 {n-1-\vert \al \vert}\prod_{i\in I}\pd {a_i} {\al_i},\\
\endcases
$$
where $\pd 1 {n-\vert \al \vert}\prod_{i\in I}\pd {a_i} {\al_i}\in
\gna n {F_S}$.
\smallskip
Let $\delta_n:A_S(n)@>>>A_S(n-1)$ be the natural projection and denote
again by $\delta_n$ its restriction to $C_S(n)$. The following diagram
is showed to be a commutative diagram in the category of graded rings:
$$
\CD \gna n {F_S}@>\rho_n>> \gna {n-1} {F_S}\\
@V\pi_nVV                      @VV\pi_{n-1}V\\
C_S(n)@>\delta_n>>C_S(n-1).\\
\endCD
$$

Take the graded inverse limit (with respect to $n$) of 
$(\gna n {F_S},\rho_n)$. This is showed to be isomorphic, via the
previous diagram, to the graded inverse limit 
(with respect to $\delta_n$) of the rings $C_{S}(n)$.

S.\,Donkin showed that the latter is the free commutative ring 
$\bigotimes_{\mu \in \psi} \Lambda_{\mu}$, where each $\Lambda_{\mu}$ 
is a copy of the ring of symmetric functions and 
$\psi$ is the set of (equivalence classes with
respect to cyclic permutations of) primitive monomials 
(see {\cite{3}},{\cite{4}}).

Let $\sigma_n:\bigotimes_{\mu \in \psi} 
\Lambda_{\mu} @>>> \gna n {F_{S}}$ be the canonical projection from
$\bigotimes_{\mu \in \psi}\Lambda_{\mu}$ considered as the graded
inverse limit of the $\gna n {F_S}$.

The proof goes on showing that the following sequence is exact:
$$0 @>>> \langle \{e_{n+1+k}(f)\, : \, k\in \NN {\text{ and }} f\in 
F_{S}^{+} \} \rangle @>>> \bigotimes_{\mu \in \psi} 
\Lambda_{\mu} @>\sigma_n>> \gna n {F_{S}}@>>>0.$$ 
Note that it makes sense to write $e_i(f)$ since this can be expressed
as a polynomial in the $e_j(m)$, with $m$ a monomial, by means of
Amitsur's Formula (see {cite{1}}).
 
In {\cite{10}} the following is proved.

\proclaim{Theorem (Zubkov)}
The following sequence is exact:
$$0 @>>> \langle \{e_{n+1+k}(f)\, : \, k\in \NN {\text{ and }} f\in 
F_{S}^{+} \} \rangle @>>> \bigotimes_{\mu \in \psi} 
\Lambda_{\mu} @>\theta_n>> C_S(n)@>>>0,$$
where $\theta_n$ is the canonical projection from 
$\bigotimes_{\mu \in \psi} \Lambda_{\mu}$ considered as the graded
inverse limit of the $C_S(n)$. 
\endproclaim

The result then follows by comparing the presentation of $\gna n
{F_S}$ with the one of $C_S(n)$ by means of $\pi_n$ and
its limit.

\smallskip
\proclaim{Corollary 2.2.3}
Via Theorem 1.2.1, the ring isomorphism $\gna n {F_{S}}\cong C_{S}(n)$
corresponds to the polynomial law $det \cdot j_{n} : F_{S} @>>> C_{S}(n)$.
\endproclaim
\demo{Proof}
It follows from Th.2.2.2 and Remark 2.2.1.
\endemo

\vskip20pt
\head \S 3 Proof of the Main Result \endhead
\vskip10pt
Let $p:F_{S} @>>> B$ be a homogeneous multiplicative polynomial 
law of degree $n$, then, by 
Theorem 1.2.1 there 
exists a unique ring homomorphism $\phi:\gna n {F_{S}} @>>> B$ 
such that $p_{_{L}}=\big(1_{_{L}}\otimes \phi\big)\cdot \gamma_{n,L}$ 
for any commutative ring $L$. But 
$\gna n {F_{S}}\cong C_{S}(n)$ by Theorem 2.2.2 
and Corollary 2.2.3 gives the arrows making the following diagram to 
commute:
$$
\CD  F_S @>p>> B \\
@V j_n VV        @AA\phi A \\
j_n(F_S) @> \det >> C_S(n).
\endCD
$$

\enddocument

\bye